\documentclass[10pt]{article}
\usepackage{amssymb, amsmath, url, graphicx, setspace, geometry}
\usepackage{calrsfs}
\usepackage{wasysym}
\usepackage{color}
\def\3{\subset }
\def\4{\subseteq }
\def\<{\left<}
\def\>{\right>}

\def\bit{\begin{itemize}}
\def\eit{\end{itemize}}
\def\3{\subset }
\def\4{\subseteq }

\def\0{\leqno}

\def\barr{\begin{array}}
\def\earr{\end{array}}

\def\Z{{\rlap{$\kern2pt{\rm Z}$}{\rm Z}\,}}
\def\bld#1#2{{\buildrel{#1}\over{#2}}}
\def\st#1#2{{\mathrel{\mathop{#2}\limits_{#1}}{}\!}}
\def\stb#1#2#3{{\st{{#1}}{\bld{{#2}}{#3}}{}\!}}
\def\xmare#1#2{\stb{#1}{#2}{\mbox{\Huge$\times$}}}

\title{\bf Finite groups with two relative subgroup commutativity degrees}
\author{Mihai-Silviu Lazorec and Marius T\u arn\u auceanu}
\date{January 27, 2018}

\begin{document}

\maketitle

\begin{abstract}
In this paper we show that there is an infinite number of finite groups with two relative subgroup commutativity degrees. Also, we indicate a sufficient condition such that a finite group has at least three relative subgroup commutativity degrees and we prove that $D_6$ is the only finite dihedral group with two relative commutativity degrees. Finally, we study the density of the set containing all subgroup commutativity degrees of finite groups.
\end{abstract}

\noindent{\bf MSC (2010):} Primary 20D60, 20P05; Secondary 20D30, 20F16, 20F18.

\noindent{\bf Key words:} subgroup commutativity degree, relative
subgroup commutativity degree, subgroup lattice.

\section{Introduction}

The large amount of results (see \cite{2}, \cite{4}-\cite{8}, \cite{10}) obtained during the last decades using the commutativity degree of a finite group $G$, determined the introduction and study of some new probabilistic aspects of finite groups theory. In this paper, we will focus on two such concepts that were introduced in \cite{15} and \cite{16}. More exactly, we will prove some results related to the relative subgroup commutativity degree of a subgroup $H$ of $G$ defined by
$$sd(H,G)=\frac{1}{|L(H)||L(G)|}|\lbrace (H_1, G_1) \in L(H)\times L(G) \ | \ H_1G_1=G_1H_1\rbrace|,$$ which is a generalization of the subgroup commutativity degree of $G$ that is given by
$$sd(G)=\frac{1}{|L(G)|^2}|\lbrace (H,K)\in L(G)^2 \ | \ HK=KH\rbrace|,$$ where $L(G)$ denotes the subgroup lattice of $G$.

We recall that if $G=\xmare{i=1}k G_i$, where $(G_i)_{i=\overline{1,k}}$ is a family of finite groups having coprime orders, the subgroup lattice of $G$ is decomposable. Hence, a subgroup $H$ of $G$, is written as $H=\xmare{i=1}k H_i$, where $H_i \in L(G_i), \ \forall \ i=1,2,\ldots,k$. Moreover, we have
$$sd(G)=\prod\limits_{i=1}^k sd(G_i) \ \text{and} \ sd(H,G)=\prod\limits_{i=1}^k sd(H_i,G_i).$$
A relevant property of the function $$f:L(G) \longrightarrow [0,1], f(H)=sd(H,G),$$
is that it is constant on each conjugacy class of subgroups of $G$. It is easy to see that $|Im \ f|=1$ if and only if $G$ is an Iwasawa group (i.e. a nilpotent modular group). Having in mind these results, our aim is to study the following class of groups
$$\mathcal{C}=\lbrace G = \text{finite group} \ | \ |Im \ f|=2\rbrace.$$

The paper is organized as follows. Section 2 describes all finite groups with two relative commutativity degrees under the additional assumption that they have one or two conjugacy classes of non-normal subgroups. A consequence of this study is the fact that $\mathcal{C}$ contains an infinite number of groups. We prove a criterion that indicates which groups do not belong to $\mathcal{C}$ in Section 3. Using this result, we detect some well known classes of groups which are not included in $\mathcal{C}$. Also, we show that $S_3\cong D_6$ is the unique finite dihedral group contained in $\cal{C}$. In Section 4, we study the density of the sets $\lbrace sd(G) \ | \ G = \text{finite group}\rbrace$ and $\lbrace sd(H,G) \ | \ G = \text{finite group}, H\in L(G)\rbrace$ in $[0,1].$ We end our paper by suggesting some open problems in the last section.

Most of our notation is standard and will usually not be repeated here. Elementary notions and results on groups can be found in \cite{3,12}. For subgroup lattice concepts we refer the reader to \cite{11,13,14}.

\section{The connection between $\mathcal{C}$ and the conjugacy classes of non-normal subgroups of a finite group}

For a finite group $G$, we will denote the number of its conjugacy classes of non-normal subgroups by $\gamma(G)$. It is obvious that $sd(\lbrace 1\rbrace,G)=1$. Therefore, it is natural to think that a group $G$ may belong to $\mathcal{C}$ if $\gamma(G)=1$ or $\gamma(G)=2$. This additional hypothesis will help us in our study because all finite groups having one or two conjugacy classes of non-normal subgroups were described in \cite{1} and \cite{9}, respectively. More exactly, the following two theorems were proved.\\

\textbf{Theorem 2.1.} \textit{Let $G$ be a finite group. Then $\gamma(G)=1$ if and only if $G$ is isomorphic to one of the following groups:
\begin{itemize}
\item[(1)] $N \rtimes P$, where $N \cong \mathbb{Z}_p, P \cong \mathbb{Z}_{q^n}, [N,\Phi(P)]=1$ and $p,q$ are primes such that $q \ | \ p-1$;
\item[(2)] $M(p^n)=\langle x, y \ | \ x^{p^{n-1}}=y^p=1, x^y=x^{1+p^{n-2}}\rangle,$ where $p$ is a prime and $n\geq 3$ if $p\geq 3$ or $n\geq 4$ if $p=2$.
\end{itemize}}

\textbf{Theorem 2.2.} \textit{Let $G$ be a finite group. Then $\gamma(G)=2$ if and only if $G$ is isomorphic to one of the following groups:
\begin{itemize}
\item[(1)] $A_4$;
\item[(2)] $\langle x,y \ | \ x^q=y^{p^n}=1, x^y=x^k\rangle$, where $p,q$ are prime numbers such that $p^2|q-1$, $n>1$ and $k^{p^2}\equiv 1 \ (mod \ q)$ with $k\ne 1$;
\item[(3)] $\langle x,y,z \ | \ x^r=y^{p^n}=z^q=[x,z]=[y,z]=1, x^y=x^k\rangle$, where $p,q,r$ are prime numbers such that $p\ne q, q\ne r, p|r-1$ and $k^{p}\equiv 1 \ (mod \ r)$ with $k\ne 1$;\,\footnote{Note that the condition $q\ne r$ must be added to the group (3) in Theorem I of \cite{9}, as shows the example $G=S_3\times\mathbb{Z}_3$; in this case we have $q=r=3$, but $\gamma(G)=3$.}
\item[(4)] $\langle x,y \ | \ x^{q^2}=y^{p^n}=1, x^y=x^k \rangle,$ where $p,q$ are primes such that $p|q-1$ and $k^p \equiv 1 \ (mod \ q^2)$ with $k\ne 1$;
\item[(5)] $M(p^n)\times \mathbb{Z}_q$, where $p,q$ are primes such that $p\ne q$ and $n\geq 3$ if $p\geq 3$ or $n\geq 4$ if $p=2$;
\item[(6)] $\mathbb{Z}_4 \rtimes \mathbb{Z}_4$;
\item[(7)] $Q_{16}$;
\item[(8)] $\langle x,y \ | \ x^4=y^{2^n}=1, y^x=y^{1+2^{n-1}}\rangle$, where $n\geq 3$;
\item[(9)] $D_8$.
\end{itemize}
}

Since the above classifications are complete, our aim is to pass through the two lists of groups and find the ones that are contained in $\mathcal{C}$. Some answers are provided by our next two results.\\

\textbf{Theorem 2.3.} \textit{Let $G$ be a finite group such that $\gamma(G)=1$. Then $G \in \mathcal{C}$ if and only if $G\cong S_3$.}

\textbf{Proof.}
We saw that there are only two classes of groups having the property $\gamma(G)=1$. The modular $p$-groups $M(p^n)$ are Iwasawa, so $|Im \ f|=1$. Consequently, they are not contained in $\mathcal{C}$. We assume that $G\cong \mathbb{Z}_p \rtimes \mathbb{Z}_{q^n}$, where $p,q$ are primes such that $q|p-1$. There are $p$ subgroups isomorphic to $\mathbb{Z}_{q^n}$ and they form the conjugacy class of non-normal subgroups of $G$. Two different conjugates cannot commute since $G$ does not contain a subgroup of order $q^{n+1}$. Therefore, $$sd(\mathbb{Z}_{q^n},G)=\frac{|L(\mathbb{Z}_{q^{n-1}})||L(G)|+|N(G)|+1}{|L(\mathbb{Z}_{q^n})||L(G)|}=\frac{n(2n+p+1)+2(n+1)}{(n+1)(2n+p+1)},$$ while
$$sd(G)=\frac{|N(G)||L(G)|+p(|N(G)|+1)}{|L(G)|^2}=\frac{(2n+1)(2n+p+1)+2p(n+1)}{(2n+p+1)^2}.$$ Then $$G\in \mathcal{C} \Leftrightarrow sd(\mathbb{Z}_{q^n},G)=sd(G) \Leftrightarrow p=2+\frac{1}{n} \ \text{or} \ p=1.$$ But $p,q$ are primes such that $q|p-1$, so $(p,q,n)=(3,2,1)$. Consequently $G\in \mathcal{C}$ if and only if $G \cong S_3$. For any other triplet $(p,q,n)$, we have $|Im \ f|=3.$
\hfill\rule{1,5mm}{1,5mm} \\

\textbf{Theorem 2.4.} \textit{Let $G$ be a finite group with $\gamma(G)=2$. Then $G\in \mathcal{C}$ if and only if $G\cong S_3 \times \mathbb{Z}_q$, where $q\geq 5$ is a prime number.}

\textbf{Proof.} We are moving our attention to the list of groups provided by \textbf{Theorem 2.2.}. First, we will inspect the nilpotent groups that have two non-normal conjugacy classes of subgroups. It is easy to see that if $G\cong Q_{16}$ or $G\cong D_8$, we obtain that $|Im \ f|=3$. Also, if $G\cong \mathbb{Z}_4\rtimes \mathbb{Z}_4$, we get $|Im \ f|=4$. Going further, a direct product of two Iwasawa groups is an Iwasawa group, so $|Im \ f|=1$ for a group $G\cong M(p^n)\times \mathbb{Z}_q$, where the triplet $(p,q,n)$ is chosen as it is indicated in \textbf{Theorem 2.2.}. If $G$ is a group of type (8), then it is a nonhamiltonian 2-group that contains an abelian normal subgroup $N\cong \mathbb{Z}_{2^n}$. According to \textbf{Theorem 9} of \cite{13}, since $\frac{G}{N}$ is cyclic and there exist $x\in G$ and an integer $m=n-1\geq 2$ such that $G=\langle N,x\rangle$ and $g^x=g^{1+2^m}, \ \forall g\in N$, we deduce that $G$ is modular. Consequently, $G$ is an Iwasawa group and we have $|Im \ f|=1$.

The groups that are not covered yet are non-nilpotent. It is easy to see that $|Im \ f|=5$ for $G\cong A_4$. If $G$ is a group of type (2), then $G\cong \mathbb{Z}_q\times \mathbb{Z}_{p^n}$, where $p,q$ are primes such that $p^2|q-1$ and $n>1$. There are $q$ subgroups isomorphic to $\mathbb{Z}_{p^n}$ and $q$ subgroups isomorphic to $\mathbb{Z}_{p^{n-1}}$ and they form the two non-normal conjugacy classes of subgroups of $G$. We will denote them by $C_1$ and $C_2$, respectively. Since $\mathbb{Z}_{p^n}$ is a cyclic $p$-group, it contains one subgroup isomorphic to $\mathbb{Z}_{p^{n-1}}$. Therefore, two different conjugates contained in $C_2$ cannot commute. Also, the same thing can be said about two different conjugates taken from $C_1$. Then, we have
$$sd(\mathbb{Z}_{p^{n-1}},G)=\frac{|L(\mathbb{Z}_{p^{n-2}})||L(G)|+|N(G)|+2}{|L(\mathbb{Z}_{p^{n-1}})||L(G)|}=\frac{(n-1)(n+q)+n+1}{n(n+q)},$$
$$sd(\mathbb{Z}_{p^n},G)=\frac{|L(\mathbb{Z}_{p^{n-2}})||L(G)|+2(|N(G)|+2)}{|L(\mathbb{Z}_{p^n})||L(G)|}=\frac{(n-1)(n+q)+2(n+1)}{(n+1)(n+q)},$$
$$sd(\mathbb{Z}_q\rtimes\mathbb{Z}_{p^{n-1}},G)=\frac{(2n-1)|L(G)|+q(|N(G)|+2)}{|L(\mathbb{Z}_q\rtimes\mathbb{Z}_{p^{n-1}})||L(G)|}=\frac{(2n-1)(n+q)+q(n+1)}{(2n+q-1)(n+q)},$$
$$sd(G)=\frac{|N(G)||L(G)|+2q(|N(G)|+2)}{|L(G)|^2}=\frac{n(n+q)+q(n+1)}{(n+q)^2}.$$
Now, to find $|Im \ f|$, we must study when two of the above (relative) subgroup commutativity degrees are equal, under the assumptions that $p,q$ are primes such that $p^2|q-1$ and $n>1$. After making some computation, we deduce that $|Im \ f|=4$ if and only if $(p,q,n)=(2,5,4)$ and $|Im \ f|=5$ if and only if $(p,q,n)\ne (2,5,4).$

A group $G$ of type (4) is isomorphic to $\mathbb{Z}_{q^2}\rtimes \mathbb{Z}_{p^n}$, where $p,q$ are primes such that $p|q-1$ and $n$ is a pozitive integer. There are $q$ subgroups isomorphic to $\mathbb{Z}_q \rtimes \mathbb{Z}_{p^n}$ which form a non-normal conjugacy class of subgroups. Each of them have $q$ subgroups isomorphic to $\mathbb{Z}_{p^n}$, so the second non-normal conjugacy class contains $q^2$ subgroups. The reader can check if two non-normal subgroups commute using similar techniques with the ones applied previously. We obtain the following results:
$$sd(\mathbb{Z}_{p^n},G)=\frac{|L(\mathbb{Z}_{p^{n-1}})||L(G)|+|N(G)|+2}{|L(\mathbb{Z}_{p^n})||L(G)|}=
\frac{n(3n+q^2+q+1)+3(n+1)}{(n+1)(3n+q^2+q+1)},$$
\begin{align*}
sd(\mathbb{Z}_q\rtimes \mathbb{Z}_{p^n},G) &=\frac{|L(\mathbb{Z}_{p^{n-1}q})||L(G)|+q(|N(G)|+2)+|N(G)|+q+1}{|L(\mathbb{Z}_q\rtimes \mathbb{Z}_{p^n})||L(G)|}\\
&= \frac{2n(3n+q^2+q+1)+3q(n+1)+3n+q+2}{(2n+q+1)(3n+q^2+q+1)},
\end{align*}
\begin{align*}
sd(G) &=\frac{|N(G)||L(G)|+q^2(|N(G)|+2)+q(|N(G)|+q+1)}{|L(G)|^2}\\
&=\frac{(3n+1)(3n+q^2+q+1)+3q^2(n+1)+q(3n+q+2)}{(3n+q^2+q+1)^2}.
\end{align*}
The properties of the triplet $(p,q,n)$ lead to stating that the above 3 (relative) subgroup commutativity degrees are different. Consequently, $|Im \ f|=4$ for a group of type (4).

Finally, if $G$ is a group of type (3), then $G\cong (\mathbb{Z}_r \rtimes \mathbb{Z}_{p^n})\times \mathbb{Z}_q$, where $p,q,r$ are prime numbers such that $p\ne q, q\ne r, p|r-1$ and $n$ is a positive integer. We remark that the group $\mathbb{Z}_r \rtimes \mathbb{Z}_{p^n}$ is exactly the non-nilpotent group having one conjugacy class of non-normal subgroups described in \textbf{Theorem 2.1.}. Also, it is obvious that $(rp,q)=1$. Consequently, the first conjugacy class of non-normal subgroups contains $r$ subgroups isomorphic to $\mathbb{Z}_{p^n}$, while the second class is formed by $r$ subgroups isomorphic to $\mathbb{Z}_{p^n}\times \mathbb{Z}_q$. Moreover, we have
$$sd(\mathbb{Z}_{p^n},G)=sd(\mathbb{Z}_{p^n},\mathbb{Z}_r \rtimes \mathbb{Z}_{p^n})sd(\lbrace 1\rbrace, \mathbb{Z}_q)=\frac{n(2n+r+1)+2(n+1)}{(n+1)(2n+r+1)},$$
$$sd(\mathbb{Z}_{p^n}\times \mathbb{Z}_q, G)=sd(\mathbb{Z}_{p^n}, \mathbb{Z}_r\rtimes \mathbb{Z}_{p^n})sd(\mathbb{Z}_q)=\frac{n(2n+r+1)+2(n+1)}{(n+1)(2n+r+1)},$$
$$sd(\mathbb{Z}_r\rtimes \mathbb{Z}_{p^n},G)=sd(\mathbb{Z}_r\rtimes \mathbb{Z}_{p^n})sd(\lbrace 1\rbrace, \mathbb{Z}_q)=\frac{(2n+1)(2n+r+1)+2r(n+1)}{(2n+r+1)^2},$$
$$sd(G)=sd(\mathbb{Z}_r\rtimes \mathbb{Z}_{p^n})sd(\mathbb{Z}_q)=\frac{(2n+1)(2n+r+1)+2r(n+1)}{(2n+r+1)^2}.$$
Then $$G\in \mathcal{C} \Longleftrightarrow sd(G)=sd(\mathbb{Z}_{p^n},G) \Longleftrightarrow r=2+\frac{1}{n} \ \text{or} \ r=1.$$
Following the properties of $p,q,r$ and $n$, we have $G\in \mathcal{C}$ if and only if $G\cong S_3\times \mathbb{Z}_q$, where $q$ is a prime number such that $q\geq 5$. We finish our proof by stating that for a group of type (3), we get $|Im \ f|=3$ for all other possible choices of the quadruplet $(p,q,r,n)$.
\hfill\rule{1,5mm}{1,5mm} \\

Looking to our last proof, we deduce that even if $\gamma(G)=n$, where $n\geq 3$ is an integer, we always find a group $G$ that belongs to $\mathcal{C}$. An example would be $G\cong S_3\times \mathbb{Z}_{5^{n-1}}$. Also, our previous proof suggests that we can generalize our last result in the following way.\\

\textbf{Corollary 2.5.} \textit{Let $G$ be a finite Iwasawa group such that $(6, |G|)=1$. Then $S_3 \times G$ is contained in $\mathcal{C}$.}

\section{How can we find finite groups that do not belong to $\mathcal{C}$?}

In the previous section, we used the number of conjugacy classes of non-normal subgroups of a finite group $G$ to establish that $\mathcal{C}$ contains an infinite number of groups. Our next aim is to find a condition which guarantees that a finite group $G$ is not a part of $\mathcal{C}$. We will denote by $N(G)$ the normal subgroup lattice of $G$ and, for a subgroup $H$ of $G$, we will consider the following set $$C(H)=\lbrace K\in L(G) \ | \ HK=KH\rbrace.$$

\textbf{Proposition 3.1.} \textit{Let $G$ be a finite group. If
$sd(G)<\frac{1}{2}+\frac{|N(G)|+1}{2|L(G)|},$ then $|Im \ f|>2.$}

\textbf{Proof.} Assume that $|Im \ f|\leq 2$. Then $|Im \ f|=2$ since $|Im \ f|=1$ would imply that $G$ is Iwasawa and this is contradicting our hypothesis. Let $H$ be a minimal subgroup of $G$ such that $sd(H,G)\ne 1$. Then $sd(H,G)=sd(G)$ as a consequence of the fact that $|Im \ f|=2$. On the other hand, we have $sd(K,G)=1$ for any proper subgroup $K$ of $H$. Hence,
$$sd(G)=sd(H,G) = \frac{1}{|L(H)||L(G)|}\sum\limits_{K\in L(H)}|C(K)|=\frac{1}{|L(H)||L(G)|}[(|L(H)|-1)|L(G)|+|C(H)|]$$
\begin{align*}
 \ \ \ \ \ \ \ \ \ \ \ \ \ \ \ \ \ \ \ \ \ \ \ &\geq \frac{1}{|L(H)||L(G)|}[(|L(H)|-1)|L(G)|+|N(G)|+1]=1-\frac{|L(G)|-|N(G)|-1}{|L(H)||L(G)|}\\
&\geq 1-\frac{|L(G)|-|N(G)|-1}{2|L(G)|}=\frac{1}{2}+\frac{|N(G)|+1}{2|L(G)|},
\end{align*}
which again contradicts our hypothesis. This argument completes our proof.
\hfill\rule{1,5mm}{1,5mm} \\

We remark that the upper bound $\frac{1}{2}+\frac{|N(G)|+1}{2|L(G)|}$ is the best possible one since, if $G\cong S_3$, we have $sd(G)=\frac{1}{2}+\frac{|N(G)|+1}{2|L(G)|}$, but $|Im \ f|=2$ as we already saw. Also, our previous result mainly states that a group $G$ may be contained in $\mathcal{C}$ if its subgroup commutativity degree is "large". Consequently, we expect that a class of groups is not included in $\mathcal{C}$ if its subgroup commutativity degree vanishes asymptotically. Some examples of well known classes of groups having this property were indicated in \cite{15}. They are\\
-the dihedral groups $$D_{2^n}=\langle x,y \ | \ x^{2^{n-1}}=y^2=1, yxy^{-1}=x^{2^{n-1}-1}\rangle, n\geq 3,$$
-the generalized quaternion groups $$Q_{2^n}=\langle x,y \ | \ x^{2^{n-1}}=y^4=1, yxy^{-1}=x^{2^{n-1}-1}\rangle, n\geq 3,$$
-the quasi-dihedral groups $$S_{2^n}=\langle x,y \ | \ x^{2^{n-1}}=y^2=1, y^{-1}xy=x^{2^{n-2}-1}\rangle, n\geq 4.$$

\textbf{Corollary 3.2.} \textit{The sets $\lbrace D_{2^n}\rbrace_{n\geq 3}$, $\lbrace Q_{2^n}\rbrace_{n\geq 3}$ and $\lbrace S_{2^n}\rbrace_{n\geq 4}$ are not included in $\mathcal{C}$.}

\textbf{Proof.} We will check if the inequality
\begin{align}\label{r1}
sd(G)<\frac{1}{2}+\frac{|N(G)|+1}{2|L(G)|}
\end{align}
holds for a finite group $G$ contained in any of the mentioned classes of groups. The following explicit formulas were indicated in \cite{15} and \cite{17}:
$$sd(D_{2^n})=\frac{(n-2)2^{n+2}+n2^{n+1}+(n-1)^2+8}{(n-1+2^n)^2},$$
$$sd(Q_{2^n})=\frac{(n-3)2^{n+1}+n2^n+(n-1)^2+8}{(n-1+2^{n-1})^2},$$
$$sd(S_{2^n})=\frac{(n-3)2^{n+1}+n2^n+(3n-2)2^{n-1}+(n-1)^2+8}{(n-1+3\cdot 2^{n-2})^2},$$
$$|N(D_{2^n})|=|N(Q_{2^n})|=|N(S_{2^n})|=n+3.$$
Using a mathematical software, we deduce that inequality (\ref{r1}) is satisfied for a finite group $G$ contained in any of the sets $\lbrace D_{2^n}\rbrace_{n\geq 5}$, $\lbrace Q_{2^n}\rbrace_{n\geq 6}$ and $\lbrace S_{2^n}\rbrace_{n\geq 5}$. Hence, for such a group we have $|Im \ f|>2$ as \textbf{Proposition 3.1.} indicates. Consequently, $G$ is not an element of $\mathcal{C}$. Moreover, it is easy to see that $$|Im \ f|=\begin{cases} 1 &\mbox{, \ if \ } G\cong Q_8 \\ 3 &\mbox{, \ if \ } G\cong D_8, G\cong Q_{16}, G\cong S_{16} \\ 4 &\mbox{, \ if \ } G\cong D_{16}, G\cong Q_{32} \end{cases},$$ this being a result which completes our proof.
\hfill\rule{1,5mm}{1,5mm} \\

\textbf{Theorem 3.3.}\textit{$S_3\cong D_6$ is the unique finite dihedral group contained in $\cal C$.}

\textbf{Proof.}
   Assume that the dihedral group $D_{2n}=\langle x,y\mid x^n=y^2=1, yxy=x^{-1}\rangle$ belongs to $\cal C$. Then $sd(H,D_{2n})=sd(K,D_{2n})$ for every two non-Iwasawa subgroups $H$ and $K$ of $D_{2n}$. The subgroup structure of $D_{2n}$ is well-known: given a divisor $r$ or $n$, $D_{2n}$ possesses a subgroup isomorphic to $\mathbb{Z}_r$, namely $H^r_0=\langle x^{\frac{n}{r}}\rangle$, and $\frac{n}{r}$ subgroups isomorphic to $D_{2r}$, namely $H^r_i=\langle x^{\frac{n}{r}},x^{i-1}y\rangle$, $i=1,2,...,\frac{n}{r}$\,. Then $|L(D_{2n})|=\tau(n)+\sigma(n)$, where $\tau(n)$ and $\sigma(n)$ are the number and the sum of all divisors of $n$, respectively. Also, by \cite{15} we have $|C(H^r_i)|=\tau(n)+x^r_i$, where $x^r_i$ is the number of solutions of
\begin{align}
\frac{n}{[r,s]}\mid 2(i-j) \mbox{ with } s\!\mid\! n \mbox{ and } j=1,2,...,\frac{n}{s}.
\end{align}
Since $n$ cannot be of type $2^m$ (see \textbf{Corollary 3.2.}), we distinguish the following two cases.

{\bf Case 1.} $n$ is odd, say $n=p_1^{\alpha_1}p_2^{\alpha_2}\cdots p_k^{\alpha_k}$ with each $p_i$ an odd prime.\\
   By taking $H=H^1_1$ and $K=H^{p_i}_1$, one obtains $$|C(H)|=2\tau(n) \mbox{ and } |C(K)|=2\tau(n)+(p_i-1)\tau(\frac{n}{p_i^{\alpha_i}}),$$implying that
   $$sd(H,D_{2n})=\frac{3\tau(n)+\sigma(n)}{2(\tau(n)+\sigma(n))} \mbox{ and } sd(K,D_{2n})=\frac{(2p_i+4)\tau(n)+2\sigma(n)+(p_i-1)\tau(\frac{n}{p_i^{\alpha_i}})}{(p_i+3)(\tau(n)+\sigma(n))}\,.$$
   These equalities lead to
\begin{align}
  \sigma(n)-\tau(n)=2\tau(\frac{n}{p_i^{\alpha_i}}).
  \end{align}
Since $i$ is arbitrary, we infer that $\alpha_1=\alpha_2=...=\alpha_k=\alpha$, and (3) becomes
\begin{align*}
\prod_{i=1}^k\frac{p_i^{\alpha+1}-1}{p_i-1}=(\alpha+1)^{k-1}(\alpha+3).
\end{align*}
It is easy to see that $$\frac{p_i^{\alpha+1}-1}{p_i-1}\geq(\alpha+1)^2,\, \forall\, i=1,2,...,k,$$which implies that $$(\alpha+1)^{k-1}(\alpha+3)\geq(\alpha+1)^{2k}, \ \text{i.e.} \ (\alpha+3)\geq(\alpha+1)^{k+1}.$$From this inequality it follows that $k=\alpha=1$ and $p_1=3$. Thus $n=3$, as desired.

   {\bf Case 2.} $n=2^mn'$, where $m\geq 1$ and $n'=p_1^{\alpha_1}p_2^{\alpha_2}\cdots p_k^{\alpha_k}$ is odd.\\
   Similarly, by taking $H=H^1_1$ and $K=H^{p_i}_1$, one obtains $$|C(H)|{=}\tau(n)+(2m+1)\tau(n') \mbox{ and } |C(K)|=\tau(n)+(2m+1)\!\left[\tau(n')+(p_i-1)\tau(\frac{n'}{p_i^{\alpha_i}})\right]\!,$$implying that
   $$sd(H,D_{2n})=\frac{2\tau(n)+\sigma(n)+(2m+1)\tau(n')}{2(\tau(n)+\sigma(n))}$$and  $$sd(K,D_{2n})=\frac{(p_i+3)\tau(n)+2\sigma(n)+(2m+1)(p_i+1)\tau(n')+(2m+1)(p_i-1)\tau(\frac{n'}{p_i^{\alpha_i}})}{(p_i+3)(\tau(n)+\sigma(n))}\,,$$respectively. Then
\begin{align}
\sigma(n)-(2m+1)\tau(n')=2(2m+1)\tau(\frac{n'}{p_i^{\alpha_i}}),
\end{align}
 and again we have $\alpha_1=\alpha_2=...=\alpha_k=\alpha$. Thus (4) becomes
\begin{align}
(2^{m+1}-1)\sigma(n')=(2m+1)(\alpha+1)^{k-1}(\alpha+3).
\end{align}
But we have $2^{m+1}-1\geq 2m+1$ with equality if and only if $m=1$, and $\sigma(n')\geq(\alpha+1)^{k-1}(\alpha+3)$ with equality if and only if $k=\alpha=1$ and $p_1=3$. So, (5) leads to $n=6$. In this case we observe that $$sd(H,D_{12})=sd(K,D_{12})=\frac{13}{16}\neq\frac{101}{128}=sd(D_{12},D_{12})=sd(D_{12})$$(i.e. the function $f$ has at least three distinct values $1$, $\frac{13}{16}$, and $\frac{101}{128}$), a contradiction.

Hence $n=3$, completing the proof.
\hfill\rule{1,5mm}{1,5mm}

\section{Density results for (relative) subgroup commutativity degree}

Let $\alpha\in[0,1]$. We will consider the following two sets
$$A=\lbrace sd(G) \ | \ G = \text{finite group}\rbrace \ \text{and} \ B=\lbrace sd(H,G) \ | \ G = \text{finite group}, H\in L(G)\rbrace.$$
An interesting question that appeared after the subgroup commutativity degree of a finite group $G$ was introduced, was if there exists a sequence of groups $(G_n)_{n\in \mathbb{N}}$, such that $\displaystyle\lim_{n\to \infty} sd(G_n)=\alpha$. In other words, is the set $A$ dense in $[0,1]$? We will show that the answer is negative and an important role in our proof will be played by the following preliminary result.\\

\textbf{Proposition 4.1.} \textit{The set $B$ is dense in $[0,1]$.}

\textbf{Proof.} Let $\alpha\in [0,1]$. It is obvious that $\alpha=1$ is a limit point. If $\alpha=0$, then we choose the sequence $(D_{2^n},D_{2^n})_{n\geq 3}$ and we have $\displaystyle\lim_{n\to \infty} sd(D_{2^n},D_{2^n})=\displaystyle\lim_{n\to \infty} sd(D_{2^n})=0,$ as it was proved in \cite{15}. Further, we consider $\alpha=\frac{a}{b}\in \mathbb{Q}\cap (0,1)$, where $a,b$ are positive integers with $a<b$. The proof of \textbf{Theorem 2.3.} provides the explicit formula that allows us to compute
$sd(\mathbb{Z}_{q^n},\mathbb{Z}_p\rtimes \mathbb{Z}_{q^n})$, where $p,q$ are primes such that $q|p-1$ and $n$ is a positive integer. We remark that $$\displaystyle\lim_{p\to \infty}sd(\mathbb{Z}_{q^n},\mathbb{Z}_p\rtimes \mathbb{Z}_{q^n})=\displaystyle\lim_{p\to \infty}\frac{n(2n+p+1)+2(n+1)}{(n+1)(2n+p+1)}=\frac{n}{n+1}.$$
Further, we consider the sequence $(q_i)_{i\in\mathbb{N}}$, where each $q_i$ is a prime number of the form $4k+3$, with $k\in\mathbb{N}$. Since $(4q_i,1)=1$, there is a sequence of primes $(p_i)_{i\in\mathbb{N}}$ of the form $4hq_i+1$. In this way, for each prime $q_i$, we find a prime $p_i$ such that $q_i|p_i-1$. Moreover the sequences $(p_i)_{n\in\mathbb{N}}, (q_i)_{n\in\mathbb{N}}$ are strictly increasing and $p_i\ne q_j, \forall \ i, j\in \mathbb{N}$. Now, let $(k_n^1), (k_n^2), \ldots, (k_n^{b-a})$ some strictly increasing and disjoint subsequences of $\mathbb{N}$. From our previous discussion about prime numbers, we can choose the sequences $(H_j^n,G_j^n)_{n\in\mathbb{N}}=(\mathbb{Z}_{q_{k_n^j}^{a+j-1}},\mathbb{Z}_{q_{k_n^j}^{a+j-1}}\rtimes \mathbb{Z}_{p_{k_n^j}})_{n\in\mathbb{N}}$, where $j=1,2,\ldots, b-a$. Then $$\displaystyle\lim_{n\to\infty}sd(H_j^n,G_j^n)=\frac{a+j-1}{a+j}, \ \forall \ j=1,2,\ldots, b-a.$$
Moreover, we can build the sequence $\bigg(\xmare{j=1}{b-a} H_j^n, \xmare{j=1}{b-a} G_j^n\bigg)_{n\in\mathbb{N}}$, and use the fact that $L\bigg(\xmare{j=1}{b-a} G_j^n\bigg)$ is decomposable for each positive integer $n$, which, again, is a consequence of the above discussion about prime numbers. Finally, we get that
$$\displaystyle\lim_{n\to\infty}sd\bigg(\xmare{j=1}{b-a} H_j^n, \xmare{j=1}{b-a} G_j^n\bigg)=\prod\limits_{j=1}^{b-a}\displaystyle\lim_{n\to\infty}sd(H_j^n,G_j^n)=\prod\limits_{j=1}^{b-a}\frac{a+j-1}{a+j}=\frac{a}{b}=\alpha.$$
This implies that $(0,1)\cap\mathbb{Q}\subseteq\overline{B}$. Since $0$ and $1$ were other limit points and $B\subseteq [0,1]$, we have $[0,1]\cap \mathbb{Q}\subseteq \overline{B}\subseteq [0,1]$. Finally, we deduce that $B$ is dense in $[0,1]$ as a consequence of the fact that the closure of $[0,1]\cap \mathbb{Q}$ is $[0,1]$.
\hfill\rule{1,5mm}{1,5mm} \\

\textbf{Theorem 4.2.} \textit{The set $A$ is not dense in $[0,1]$.}

\textbf{Proof.} It is clear that $A\subseteq B$ and we proved that $\overline{B}=[0,1]$. If $B=\overline{A}$, then $B$ is closed, so $B=\overline{B}$. But $B$ is dense in $[0,1]$ and this leads us to $B=[0,1]$, which is false since any relative subgroup commutativity degree of a subgroup $H$ of $G$ can only be a rational number according to the definition of this concept. Hence, $B\ne\overline{A}$. Then there is a non-empty open set $D\subset B$ such that $D\cap A=\emptyset$. Since $D$ is an open set of $B$, it is written as $D=B\cap D_1$, where $D_1$ is an open set of $\mathbb{R}$. But $D\subset [0,1]$, so
$D=B\cap (D_1\cap [0,1])=B\cap D_2,$
where $D_2$ is an open set of $[0,1]$. We add that $D_2$ is non-empty and $D_2\cap A=\emptyset$. Otherwise, we would obtain that $D=\emptyset$ or $D\cap A\ne \emptyset$, a contradiction. Therefore, $\overline{A}\ne [0,1]$ and our proof is complete.
\hfill\rule{1,5mm}{1,5mm}

\section{Further research}

We did not fully described the set $\mathcal{C}$ and another remark is that we did not find any nilpotent group that belongs to this set. Hence, it is natural to indicate the following two open problems:\\

\textbf{Problem 5.1.} Show that $\mathcal{C}=\lbrace S_3\times G \ | \ G = \text{finite Iwasawa group such that} \ (6,|G|)=1\rbrace$.\\

\textbf{Problem 5.2.} Does $\mathcal{C}$ contains any nilpotent group?

\vspace*{3ex}
\small

\begin{minipage}[t]{7cm}
Mihai-Silviu Lazorec \\
Faculty of  Mathematics \\
"Al.I. Cuza" University \\
Ia\c si, Romania \\
e-mail: {\tt mihai.lazorec@student.uaic.ro}
\end{minipage}
\hfill
\begin{minipage}[t]{7cm}
Marius T\u arn\u auceanu \\
Faculty of  Mathematics \\
"Al.I. Cuza" University \\
Ia\c si, Romania \\
e-mail: {\tt tarnauc@uaic.ro}
\end{minipage}

\end{document}